\documentclass[11pt]{amsart}
\usepackage{graphicx}
\usepackage{amssymb, amsmath}
\newtheorem{theorem}{Theorem}[section]

\newtheorem{lemma}[theorem]{Lemma}
\newtheorem{proposition}[theorem]{Proposition}
\theoremstyle{definition}

\theoremstyle{remark}

\begin{document}
\title{Representations of Finite Polyadic Groups}
\author{\sc M. Shahryari}
\thanks{{\scriptsize
\hskip -0.4 true cm MSC(2010): 20N15
}}

\address{ Department of Pure Mathematics,  Faculty of Mathematical
Sciences, University of Tabriz, Tabriz, Iran }
\email{mshahryari@tabrizu.ac.ir}
\date{\today}

\begin{abstract}
We prove that there is a one-one correspondence between sets of
irreducible representations of a polyadic group and its Post's
cover. Using this correspondence, we generalize some well-known
properties of irreducible characters in finite groups to the case of
polyadic groups.
\end{abstract}

\maketitle

\section{Introduction}
A non-empty set $G$ together with an $n$-ary operation $f:G^n\to G$
is called an {\it $n$-ary group} or a {\it polyadic group}, if the operation $f$ is associative and for all
$x_0,x_1,\ldots,x_n\in G$ and fixed $i\in\{1,\ldots,n\}$ there
exists a unique  element $z\in G$ such that
$$
f(x_1^{i-1},z,x_{i+1}^n)=x_0.
$$
In the binary case (i.e., for $n=2$), a polyadic group is just usual
group.

In this paper, a sequence of elements $x_i,x_{i+1},\ldots,x_j$ is
denoted by $x_i^j$.  If
$x_{i+1}=x_{i+2}=\ldots=x_{i+t}=x$, then instead of $x_{i+1}^{i+t}$
we write $\stackrel{(t)}{x}$. In this convention $f(x_1,\ldots,x_n)=
f(x_1^n)$ and hence the associativity of $f$ can be formulated as
$$
f(x_1^{i-1},f(x_i^{n+i-1}),x_{n+i}^{2n-1})=
f(x_1^{j-1},f(x_j^{n+j-1}),x_{n+j}^{2n-1}),
$$
where $1\leq i, j\leq n$, and $x_1,\ldots,x_{2n-1}\in G$.

The idea of investigations of such polyadic group  goes  back to E.
Kasner's lecture \cite{Kas} at the fifty-third annual meeting of the
American Association for the Advancement of Science in 1904. But the
first paper concerning the theory of $n$-ary groups was written
(under inspiration of Emmy Noether) by W. D\"ornte in 1928 (see
\cite{Dor}). In this paper D\"ornte observed that any $n$-ary system
$(G,f)$ of the form $\,f(x_1^n)=x_1\circ x_2\circ\ldots\circ
x_n\circ b$, where $(G,\circ)$ is a group and $b$ is its fixed
element belonging to the center of $(G,\circ)$, is an $n$-ary group.
Such $n$-ary group is called {\it $b$-derived} from the group
$(G,\circ)$, and we will denote it by $der^n_b(G, \circ)$. In the
case when $b$ is the identity of $(G,\circ)$, we say that such
$n$-ary group is {\it reducible} to the group $(G,\circ )$ or {\it
derived} from $(G,\circ)$ and we denote it by $der^n(G, \circ)$. For every $n>2$ there are $n$-ary
groups which are not derived from any group. An $n$-ary group
$(G,f)$ is derived from some group if and only if it contains an
element $e$ (called an {\it $n$-ary identity}) such that
$$
 f(\stackrel{(i-1)}{e},x,\stackrel{(n-i)}{e})=x
$$
holds for all $x\in G$ and for all $i=1,\ldots,n$.

From the definition of an $n$-ary group $(G,f)$ we can directly see
that for every $x\in G$ there exists only one $z\in G$ satisfying
the equation
$$
f(\stackrel{(n-1)}{x},z)=x .
$$
This element is called {\it skew} to $x$ and is denoted by
$\overline{x}$. In a ternary group ($n=3$) derived from a binary
group $(G,\circ)$, the skew element coincides with the inverse
element in $(G,\circ)$. Thus, in some sense, the skew element is a
generalization of the inverse element in binary groups.

Nevertheless, the concept of skew elements plays a crucial role in
the theory of $n$-ary groups. Namely, as D\"ornte proved, the following theorem is true.

\begin{theorem}
In any $n$-ary group $(G,f)$ the following identities
$$
f(\stackrel{(i-2)}{x},\overline{x},\stackrel{(n-i)}{x},y)=
f(y,\stackrel{(n-j)}{x},\overline{x},\stackrel{(j-2)}{x})=y,
$$
$$
f(\stackrel{(k-1)}{x},\overline{x},\stackrel{(n-k)}{x})=x
$$
 hold for all $\,x,y\in G$, $\,2\leqslant i,j\leqslant
n$ and $1\leqslant k\leqslant n$.
\end{theorem}

Suppose $(G, f)$ is an $n$-ary group. A map $\Lambda:G \rightarrow
GL_m(\mathbb{C})$ with the property
$$
\Lambda(f(x_1^n))=\Lambda(x_1)\Lambda(x_2)\ldots\Lambda(x_n)
$$
is a {\it representation} of $G$. The function
$$
\chi(x)=Tr\ \Lambda(x)
$$
is called the corresponding {\it character} of $\Lambda$. The number
$m$ is the degree of representation. Note that, $\Lambda$ is a
representation of $(G, f)$, iff it is an $n$-ary homomorphism $G\to
der^n(GL_m(\mathbb{C}))$.

In \cite{Shah}, a joint paper with W. Dudek, we studied
representation theory of polyadic group, but representations we
dealt with in that paper were considered to have non-empty kernels.
In this paper, we study representations of polyadic groups without
that assumption, i.e. the representations we deal with in this
paper, may have empty kernels, as well. We will prove that there is
a one-one correspondence between the sets of irreducible
representations of polyadic groups and their {\it Post cover}. Using
this correspondence, we will generalize some well-known properties
of irreducible characters of finite groups to finite polyadic
groups.

\section{Generalities}
Suppose $(G,f)$ is an $n$-ary group and $a\in G$ is any fixed
element. Let
$$
G^{\ast}_a=\{ (x, i):\ x\in G, \ i\in \mathbb{Z}_{n-1}\}.
$$
We define a binary operation of $G^{\ast}_a$ by
$$
(x, i)\cdot(y, j)=(f_{\ast}(x,\stackrel{(i)}{a},y,
\stackrel{(j)}{a}, \bar{a}, \stackrel{(n-2-i\ast j)}{a}), i\ast j),
$$
where $i\ast j \equiv i+j+1\ (mod\ n-1)$, and $f_{\ast}$ indicates
that $f$ is used one or twice, depending on the value of $i\ast j$.
The set $G^{\ast}_a$ together with this operation is an ordinary
group (see \cite{Mich}), and we call it {\it Post's cover} of $(G,
f)$. The element $(\overline{a},n-2)$ is the identity of the group
$G_a^{\ast}$. The inverse element has the form
$$
( x,i)^{-1}= (
f_{\ast}(\overline{a},\stackrel{(n-2-i)}{a},\overline{x},\stackrel{(n-3)}{x},\overline{a},\stackrel{(i+1)}{a}),
k),
$$
where $k=(n-3-i)({\rm mod}\,(n-1))$. The following theorem is known
as {\it Post's coset theorem}, and it is proved in \cite{Mich} and
also in\cite{Post}.

\begin{theorem}
Suppose
$$
H=\{ (x, n-2):\ x\in G\}.
$$
Then $H\unlhd G^{\ast}_a$ and $G^{\ast}_a/H\cong \mathbb{Z}_{n-1}$.
Further, we can identify $G$ with the subset
$$
\{ (x, 0):\ x\in G\},
$$
and in under this identification, $G$ is a coset of $H$ and we have
$f(x_1^n)=x_1x_2\cdots x_n$.
\end{theorem}

It is proved that (see \cite{Mich}), Post's cover $G^{\ast}_a$ does
not depend on $a$, i.e. for any $a, b\in G$, we have
$G^{\ast}_a\cong G^{\ast}_b$.

\begin{proposition}
Suppose $A$ is an ordinary group and $a\in A$. Then for any $n\geq
2$, we have
$$
(der^n(A))^{\ast}_a\cong A\times \mathbb{Z}_{n-1}.
$$
\end{proposition}

{\bf Proof}. It is enough to suppose $a=e$, the identity element
of $A$. We have
$$
(der^n(A))^{\ast}_e=\{ (x, i):\ x\in A, i\in \mathbb{Z}_{n-1}\},
$$
and also
\begin{eqnarray*}
(x, i)\cdot (y, j)&=&(xy, i\ast j)\\
                  &=&(xy, i+j+1)\\
                  &=&(x, i)(y,j)(e,1).
\end{eqnarray*}
This shows that
$$
(der^n(A))^{\ast}_e=der^2_{(e,1)}(A\times \mathbb{Z}_{n-1}).
$$
Now, define a map $\varphi:der^2_{(e,1)}(A\times
\mathbb{Z}_{n-1})\to A\times \mathbb{Z}_{n-1}$, by $\varphi(x,
i)=(x, i+1)$. It is easy
to check that $\varphi$ is an isomorphism.\\

As a result, we see that if $a\in GL_m(\mathbb{C})$, then
$$
(der^n(GL_m(\mathbb{C})))^{\ast}_a\cong GL_m(\mathbb{C})\times
\mathbb{Z}_{n-1}.
$$

Now, let $(G, f)$ be an $n$-ary group and suppose $\Lambda:G\to
der^n(GL_m(\mathbb{C}))$ is any representation. Let $a\in G$ be
fixed. We define a new map
$$
\Lambda^{\ast}_a: G^{\ast}_a\to
(der^n(GL_m(\mathbb{C})))^{\ast}_{\Lambda(a)}
$$
by the rule
$$
\Lambda^{\ast}_a(x, i)=(\Lambda(x), i).
$$

\begin{lemma}
$\Lambda^{\ast}_a$ is a group homomorphism.
\end{lemma}

{\bf Proof}. Let $B=\Lambda(a)$. Note that we have
$\Lambda(\bar{a})=\Lambda(a)^{2-n}=\bar{B}$. Now, for any $x, y\in
G$ and $i,j \in \mathbb{Z}_{n-1}$, we have
\begin{eqnarray*}
\Lambda^{\ast}_a((x,i)\cdot(y,j))&=&\Lambda^{\ast}(f_{\ast}(x,\stackrel{(i)}{a},y,
                               \stackrel{(j)}{a}, \bar{a}, \stackrel{(n-2-i\ast j)}{a}), i\ast
                                 j)\\
&=&(\Lambda(x)\Lambda(a)^i\Lambda(y)\Lambda(a)^j\Lambda(\bar{a})\Lambda(a)^{n-2-i\ast
    j}, i\ast j)\\
&=&(\Lambda(x)\Lambda(a)^i\Lambda(y)\Lambda(a)^j\Lambda(a)^{2-n}\Lambda(a)^{n-2-i\ast
    j}, i\ast j)\\
&=&(\Lambda(x)\Lambda(a)^i\Lambda(y)\Lambda(a)^{j-i\ast
    j}, i\ast j).
\end{eqnarray*}

On the other hand,
\begin{eqnarray*}
\Lambda^{\ast}_a(x,i)\cdot\Lambda^{\ast}_a(y,j)&=&(\Lambda(x),i)\cdot
            (\Lambda(y), j)\\
          &=&(\Lambda(x)B^i\Lambda(y)B^j\bar{B}B^{n-2-i\ast j},
               i\ast j)\\
          &=&(\Lambda(x)B^i\Lambda(y)B^jB^{2-n}B^{n-2-i\ast j},
               i\ast j)\\
               &=&(\Lambda(x)B^i\Lambda(y)B^{j-i\ast j},
               i\ast j)\\
               &=&(\Lambda(x)\Lambda(a)^i\Lambda(y)\Lambda(a)^{j-i\ast
                   j}, i\ast j).
\end{eqnarray*}
This shows that $\Lambda^{\ast}_a$ is a group homomorphism.\\

Note that we have an isomorphism
$$
q:(der^n(GL_m(\mathbb{C})))^{\ast}_{\Lambda(a)}\to
(der^n(GL_m(\mathbb{C})))^{\ast}_I,
$$
where $I$ is the identity matrix. It is easy to see that
$$
q(X, i)=(X\Lambda(a)^i, i),
$$
for any $X\in GL_m(\mathbb{C})$. As we saw in the previous section,
we have also an isomorphism
$$
\varphi:(der^n(GL_m(\mathbb{C})))^{\ast}_I\to GL_m(\mathbb{C})\times
\mathbb{Z}_{n-1},
$$
such that $\varphi (X, i)=(X, i+1)$. Now, let
$\pi:GL_m(\mathbb{C})\times \mathbb{Z}_{n-1}\to GL_m(\mathbb{C})$ be
the projection. Combining all of these maps, we obtain
$$
\Lambda^{\ast}=\pi \varphi q \Lambda^{\ast}_a:G^{\ast}_a\to
GL_m(\mathbb{C}),
$$
which is an ordinary representation of $G^{\ast}_a$. Note that
\begin{eqnarray*}
\Lambda^{\ast}(x, i)&=& \pi \varphi q \Lambda^{\ast}_a(x, i)\\
                    &=&\pi \varphi q(\Lambda(x), i)\\
                    &=&\pi \varphi(\Lambda(x)\Lambda(a)^i, i)\\
                    &=&\pi (\Lambda(x)\Lambda(a)^i, i+1)\\
                    &=&\Lambda(x)\Lambda(a)^i.
\end{eqnarray*}

Conversely, suppose $\Gamma:G^{\ast}_a\to GL_m(\mathbb{C})$ is an
ordinary representation of $G^{\ast}_a$. Since $G\subseteq
G^{\ast}_a$, so by restriction we obtain a representation $\Gamma_G$
for $(G, f)$.

\begin{lemma}
The maps $\Lambda \mapsto \Lambda^{\ast}$ and $\Gamma \mapsto
\Gamma_G$ are inverse to each other.
\end{lemma}

{\bf Proof}. We have
\begin{eqnarray*}
(\Lambda^{\ast})_G(x, 0)&=&\Lambda^{\ast}(x,0)\\
                        &=&\Lambda(x)\Lambda(a)^0\\
                        &=&\Lambda(x).
\end{eqnarray*}
On the other hand
\begin{eqnarray*}
(\Gamma_G)^{\ast}(x, i)&=&\Gamma(x, 0)\Gamma(a,0)^i\\
                       &=&\Gamma(x, 0)\Gamma(a, i-1)\\
                       &=&\Gamma((x,0)\cdot(a, i-1))\\
                       &=&\Gamma(f_{\ast}(x,\stackrel{(0)}{a},a,
                               \stackrel{(i-1)}{a}, \bar{a}, \stackrel{(n-2-0\ast(i-1))}{a}),
                               0\ast (i-1))\\
                       &=&\Gamma(f_{\ast}(x,\stackrel{(i)}{a}, \bar{a},
                       \stackrel{(n-i-2)}{a}), i)\\
                       &=& \Gamma(x, i).
\end{eqnarray*}
So, the maps are inverse to each other.\\

Note that $G$ is a generating set for $G^{\ast}_a$ and hence,
$\Lambda$ is irreducible, iff $\Lambda^{\ast}$ is irreducible.
Hence, we proved;

\begin{theorem}
Let $(G,f)$ be an $n$-ary group. Then there is a one-one
correspondence between the set of all irreducible representations of
$G$ and those of $G^{\ast}_a$, for any $a\in G$. This correspondence
is the map $\Lambda \mapsto \Lambda^{\ast}$. Especially, the number
of irreducible representations of any finite polyadic group is
finite.
\end{theorem}

\section{Applications}
In this section, we apply the correspondence just we obtained, to
generalize some results of representations theory of finite groups
to the case of finite polyadic groups. In this section $(G,f)$ is a
finite $n$-ary group. Let $a\in G$ be any fixed element.

Suppose $\Lambda_1$, ..., $\Lambda_k$ are all non-equivalent
irreducible representations of $G$ with degrees
$$
d_1, d_2, \ldots, d_k.
$$
Then the set of irreducible representations of $G^{\ast}_a$ is
$$
\Lambda^{\ast}_1, \ldots, \Lambda^{\ast}_k
$$
with the same set of degrees. Since for any $i$, the order of
$G^{\ast}_a$ is divisible by $d_i$, and since we have
$$
\sum_{i=1}^kd_i^2=|G^{\ast}_a|,
$$
so we have;

\begin{theorem}
The number $(n-1)|G|$ is divisible by all $d_i$, and also we have
$$
\sum_{i=1}^kd_i^2=(n-1)|G|.
$$
\end{theorem}

Denote by $\textrm{Irr} (G,f)$, the set of all irreducible characters of $G$.
We can generalize the orthogonality property of ordinary irreducible
characters, for those elements of $\textrm{Irr}(G, f)$, which have non-empty
kernels.

\begin{theorem}
Let $\chi, \psi\in \textrm{Irr}(G,f)$ and $p\in \ker \chi$, $q\in \ker \psi$,
and $a\in G$. Then we have
$$
\frac{1}{(n-1)|G|}\sum_{i=0}^{n-2}\sum_{x\in
                     G}\chi(f(x,\stackrel{(i)}{a},
                     \stackrel{(n-i-1)}{p}))\psi(f(x,\stackrel{(i)}{a},\stackrel{(n-i-1)}{q}))^{\ast}=\delta_{\chi,
                     \psi},
$$
where $\ast$ denotes complex conjugation.
\end{theorem}

{\bf Proof}. Let $\hat{\chi}$ and $\hat{\psi}$ be the corresponding
characters of $G^{\ast}_a$.
Suppose $\Lambda$ is a representation of $G$, whose character is
$\chi$. Then
\begin{eqnarray*}
\hat{\chi}(x,i)&=&Tr\ \Lambda^{\ast}(x,i)\\
               &=&Tr\ \Lambda(x)\Lambda(a)^i\\
               &=&Tr\ \Lambda(x)\Lambda(a)^i\Lambda(p)^{n-i-1}\\
               &=&Tr\ \Lambda(f(x,\stackrel{(i)}{a},\stackrel{(n-i-1)}{p}))\\
               &=&\chi(f(x,\stackrel{(i)}{a},\stackrel{(n-i-1)}{p})).
\end{eqnarray*}
Similarly, for $\psi$, we have
$$
\hat{\psi}(x,i)=\psi(f(x,\stackrel{(i)}{a},\stackrel{(n-i-1)}{q})).
$$
Hence, we have
\begin{eqnarray*}
\delta_{\chi, \psi}&=&\delta_{\hat{\chi},\hat{\psi}}\\
                   &=&\frac{1}{|G^{\ast}_a|}\sum_{(x,i)\in G^{\ast}_a}\hat{\chi}(x,i)\hat{\psi}(x,i)^{\ast}\\
                   &=&\frac{1}{(n-1)|G|}\sum_{i=0}^{n-2}\sum_{x\in
                     G}\chi(f(x,\stackrel{(i)}{a},
                     \stackrel{(n-i-1)}{p}))\psi(f(x,\stackrel{(i)}{a},\stackrel{(n-i-1)}{q}))^{\ast}.
\end{eqnarray*}

\end{document}